\documentclass[12pt]{article}
\usepackage{graphicx}
\usepackage{mathptmx}
\usepackage{latexsym,amsmath,amssymb,amsfonts,amsthm}

\newcommand{\R}{\mathbb{R}}

\begin{document}
\setcounter{page}{1}
\title{Self-Similar Solutions to Curvature Flow of Convex Hypersurfaces}
\author{Guanghan Li$^{1\dag}$,  Isabel
Salavessa$^{2\ddag}$ and   Chuanxi Wu$^{1\dag}$}
\date{}
\protect\footnotetext{\!\!\!\!\!\!\!\!\!\!\!\!\! {\bf MSC 2000:}
Primary: 53C24 ; Secondary: 53C44.
$\dag$ Partially
supported by NSFC (No.10501011) and by NSF of Hubei Province (No.2008CDB009).\\
$\ddag$ Partially supported by FCT through the Plurianual of CFIF.}
\maketitle~~~\\[-5mm]
{\footnotesize $^1$
  School of Mathematics and Computer Science,
Hubei University, Wuhan, 430062, P. R. China,\\
{liguanghan@163.com (G.\ Li), cxwu@hubu.edu.cn (C.\ Wu)}}\\[1mm]
{\footnotesize $^2$ 
Centro de F\'{\i}sica das Interac\c{c}\~{o}es
Fundamentais, Instituto Superior T\'{e}cnico, Technical University
of Lisbon, Edif\'{\i}cio Ci\^{e}ncia, Piso 3, Av.\ Rovisco Pais,
P-1049-001 Lisboa, Portugal\\
{isabel.salavessa@ist.utl.pt}}
\begin{abstract} \noindent We classify
the self-similar solutions to a class of Weingarten curvature flow
of connected compact convex hypersurfaces, isometrically immersed
into space forms with non-positive curvature, and obtain a new
characterization of a sphere in a Euclidean space $\R^{n+1}$.
\end{abstract}
\markright{\sl \hfill G.\ Li, I.\ Salavessa and C.\ Wu \hfill}

\section{Introduction}

\renewcommand{\thesection}{\arabic{section}}
\renewcommand{\theequation}{\thesection.\arabic{equation}}
\setcounter{equation}{0}

It is a fundamental problem to classifying hypersurfaces in a
Euclidean space in classical differential geometry. For a compact
and connected hypersurface in $\R^{n+1}$, various conditions have
been obtained to guarantee that it is a standard Euclidean sphere,
and thus various characterizations of spheres have been given.

Let $X:M\rightarrow\R^{n+1}$ be a hypersurface immersed in a
Euclidean space $\R^{n+1}$. Denote by $A$ and $H$ the Weingarten
transformation and mean curvature of $M$, respectively. Assume $\bf
v$ is the unit normal vector field, then the support function of the
hypersurface $M$ is defined by
$$\mathcal{Z}=<X,{\bf v}>.$$

It is known that $M$ is a Euclidean sphere if and only if its
support function $\mathcal{Z}$ is constant and its Weingarten
transformation $A$ is not degenerate. When $M$ is oriented,
Liebmann-S\mbox{$\ddot{u}$}ss' theorem implies that, it is a
Euclidean sphere if and only if it has constant mean curvature and
its support function $\mathcal{Z}$ does not change sign. If $M$ is
closed and strictly convex, the constant mean curvature can
guarantee that it is a standard sphere. For an embedded closed
hypersurface, an interesting result of Ros \cite{r} shows that $M$
is a standard sphere if its scalar curvature is constant. For
hypersurfaces of constant Weingarten curvatures (see below for
definitions) immersed into space forms, Ecker-Huisken \cite{eh}
completely classify such hypersurfaces with non-negative sectional
curvatures.

   In 1990, Huisken \cite{h3} obtained a different characterization
of standard shperes in Euclidean spaces by studying self-similar
solutions to mean curvature flow of hypersurfaces. The mean
curvature flow is a family of evolving hypersurfaces in direction of
normal vectors, such that the evolution speed is the mean curvature.
More precisely, for a smooth oriented $n(n\ge 2)$-dimensional
manifold $M$ and $X: M\longrightarrow \R^{n+1}$ the smooth
hypersurface immersed in a Euclidean space $\R^{n+1}$, the mean
curvature flow is the following evolving problem (cf. \cite{h1})
\begin{equation}
\frac {\partial }{\partial t}X(x,t) = H(x, t){\bf v}(x,t),
                                             \quad x\in M,
\end{equation}
satisfying the initial condition $X(x,0) =X_0(x)=X(x)$, $x\in M$,
where $H(x,t)$ is the mean curvature and ${\bf v}(x,t)$ the inward
unit normal vector of $M_t=X_t(M)$ at $X(x,t)$.

It is known that, (1.1) is a contracting curvature flow, and when
the initial hypersurface is convex, Huisken \cite{h1} proved that
the solutions exist in a maximal finite time interval and converge
to a standard sphere by scaling. Later, Andrews \cite{a1} extended
this phenomenon to a class of general curvature flow, i.e.
\begin{equation}
\frac {\partial }{\partial t}X(x,t) = F(A(x, t)){\bf v}(x,t),
                                             \quad x\in M,
\end{equation}
where $F$ is a curvature function (i.e. positive and elliptic) of
homogeneous degree one of the evolving hypersurfaces satisfying
standard conditions, and $A(x,t)$ the Weingarten form of the
corresponding evolving hypersurfaces $M_t$.

If we only assume that the initial hypersurface $X_0$ has
non-negative mean curvature, the Type (I) solution to the evolving
problem (1.1) is asymptotically self-similar, i.e. the limit
hypersurface $X$ of the scaling solutions satisfies the following
equation (cf. \cite{h3})
\begin{equation}
H+<X,{\bf v}>=0.
\end{equation}
This  is a fully nonlinear elliptic equation which relates the
support function $\mathcal{Z}$ and the mean curvature $H$ of the
hypersurface $X$. Huisken \cite{h3} completely classified such
self-similar solutions to (1.3). When $M$ is compact, the only
possible case is a standard sphere, i.e. Huisken gave a new
characterization of Euclidean spheres:\\[3mm]
{\bf Proposition 1.1}~ {\it A compact and connected hypersurface
with non-negative mean curvature immersed in a Euclidean space is a
standard sphere if and only if (1.3) holds.}\\[3mm]
\indent Let $f(\lambda)$ be a function defined on a symmetric region
in $\R^n$. It is easy to see that $f$ induces a function
$F(A)=f(\lambda (A))$ defined in the set of symmetric matrices with
eigenvalues $\lambda$. When $f$ is evaluated at the vector $\lambda
(x)=\{\lambda_1(x), \cdots, \lambda_n (x)\}$, the components of
which are the principal curvatures of $M$, the hypersurface $M$ with
curvature $F$ is the so-called Weingarten hypersurface.

 If $X: M \rightarrow \R^{n+1}$ is a
sphere immersed in $ \R^{n+1}$, then there exists a constant $\tau$
such that
\begin{equation}
F+\tau <X, {\bf v}>=0,
\end{equation}
is trivially satisfied for any symmetric and homogeneous curvature
function $F$ defined in the space of positive transformations, for,
in this case the principal curvatures are all equal. We also call
hypersurfaces satisfying the condition (1.4) the self-similar
solutions of the convex curvature flow (1.2).

 Huisken's  theorem \cite{h3} says that the inverse is true for $F=H$.
In this paper, we will
consider closed immersed hypersurfaces in a space form $N^{n+1}(c)$
with non-positive curvature $c$, and show that the inverse of (1.4)
is true for a large class of Weingarten hypersurface with Weingarten
curvature $F$ satisfying some given conditions.

For this purpose, we first introduce the following functions for any
real number $c$ (cf. \cite{cgm})
\begin{eqnarray*}
sh_c(t)=\left\{\begin{array}{ll}\frac{\sin(\sqrt{c}t)}{\sqrt{c}}&\mbox{
if }c>0\\
t&\mbox{
if }c=0\\
\frac{\sinh(\sqrt{-c}t)}{\sqrt{-c}}&\mbox{ if
}c<0\end{array}\right., &~~~~\mbox{ and }~~~&
ch_c(t)=\left\{\begin{array}{ll}\cos(\sqrt{c}t)&\mbox{
if }c>0\\
1&\mbox{
if }c=0\\
\cosh(\sqrt{-c}t)&\mbox{ if }c<0\end{array}\right..
\end{eqnarray*}

Given any fixed point in the ambient space $N^{n+1}(c)$, we shall
denote by $\rho $ the distance function to the fixed point in
$N^{n+1}(c)$, and denote by $\partial _{\rho}$ the gradient of $\rho
$ in $N^{n+1}(c)$. For a hypersurface $X: M \rightarrow N^{n+1}(c)$,
let
\begin{eqnarray}
\mathcal{Z}=sh_c(\rho)<\partial _{\rho}, \bf v>,
\end{eqnarray}
where $<\cdot, \cdot>$ is the metric of the ambient space, and $\bf
v$ is again the inward unit normal vector of $M$. It is easy to see
that when $c=0$, $sh_c(\rho)\partial_{\rho}=\rho\partial_{\rho} $ is
the position vector, and therefore $\mathcal{Z}$ in (1.5) coincides
with the support function of the hypersurface.

Denote by $\Gamma_+$ the positive cone of $\R^n$, and $\Gamma(F)$ a
component of $\{\lambda :f(\lambda )\ne 0\}$ containing
$\Gamma_+$.\\[3mm]
{\bf Theorem 1.2}~ {\it Let $F(A)=f(\lambda(A))$ be a smooth
symmetric function of homogeneous degree $m\in \R/\{0\}$, and
$N^{n+1}(c)$ a Riemannian manifold of non-positive constant
curvature $c$. Suppose $X:M\rightarrow N^{n+1}(c)$  is a smooth
connected compact convex hypersurface with
 principal curvatures $\lambda=(\lambda_1, \ldots, \lambda_n)\in \Gamma_+$.
Assume the following  conditions are satisfied:\\[2mm]
 $(1)$~ On $\Gamma (F)$, $F$ is
elliptic, i.e.\ $\partial
f/\partial \lambda _i>0, \forall i=1,2 \cdots, n$.\\[1mm]
$(2)$~ One of the following holds: $(i)$ $m\ge 1$ and $f$ is convex
or concave; $(ii)$ $m<0$ and $f$ is convex or concave; $(iii)$ $n=2$
and either $m=1$, or $-7\le m<0$, or $m>1$ and $r_{\max}\le\frac12
\left(1+\sqrt{1+\frac {8}{m-1}}\right)$, where
$r=\frac{\lambda_2}{\lambda_1}\ge 1$ is the pinching ratio of the
principal curvatures, or $m< -7$ and 
$r_{\max}\le2/\left(1+\sqrt{1-\frac {8}{1-m}}\right)$.\\[2mm]
Then, if
\begin{equation}
F+\tau \mathcal{Z}=0,
\end{equation}
holds for a nonzero constant $\tau $ depending only on $n$, $X(M)$
is an umbilical sphere.}\\[3mm]
 \indent We remark that spheres are stable
solutions to contracting as well as expanding curvature flows of
convex hypersurfaces in Euclidean spaces, and therefore for $m=1$ or
$n=2$, Theorem 1.2 in fact follows from \cite{h1,a1,a2,a5}.\\[3mm]
{\bf Corollary 1.3}~ {\it  A connected compact convex hypersurface
immersed in a Euclidean space is a standard Euclidean sphere if and
only if (1.6) holds for some curvature function given in Theorem
1.2.}\\[3mm]
\indent For self-similar solutions of mean curvature flow on
arbitrary codimension, Smoczyk \cite{s} classifies such
self-shrinkers with parallel principal normal vector field. In the
case of isotropic curve flow, Andrews \cite{a3} completely
classified the homothetically shrinking solutions of (1.2), even the
curvature function is not homogeneous of degree one. For the
behavior of embedded expanding convex solutions to (1.2), there are
also complete descriptions, see \cite{a2,ct,u}, and so on.

\section{Preliminaries}
\renewcommand{\thesection}{\arabic{section}}
\renewcommand{\theequation}{\thesection.\arabic{equation}}
\setcounter{equation}{0}
\setcounter{theorem}{0}

 Let $N^{n+1}(c)$ be an $(n+1)$-dimensional space form of constant curvature
$c$, and $M$ a smooth hypersurface immersion in $N$. We will use the
same notation as in \cite{h1,a1,m3}. In particular, $\nabla $ is
the induced connection on $M$, and for a local coordinate system
$\{x^1, \cdots, x^n\}$ of $M$, $g=g_{ij}$ and $A=h_{ij}$ denote
respectively the metric and second fundamental form of $M$. Let
$g^{ij}$ denote the $(i,j)$-entry of the inverse of the matrix
$(g_{ij})$. Then $\{h_i^j\}$ where $h_i^j=h_{ik}g^{kj}$ is the
Weingarten map. The mean curvature and the squared norm of the
second fundamental form of $M$ are given by
\begin{eqnarray*}H=g^{ij}h_{ij}=h_i^i, \qquad |A|^2
=g^{ij}g^{kl}h_{ik}h_{jl}. \end{eqnarray*}  In the sequel we will
use $\lambda _i$ to denote the $i$-th principal curvature of the
hypersurface. Throughout this paper we sum over repeated indices
from $1$ to $n$ unless otherwise indicated. Raised indices indicate
contraction with the metric.

Given a symmetric smooth function $f(\lambda)$ defined in the
symmetric region of $\R^n$, the induced function $F(A)=f(\lambda
(A))$ defined in the set of symmetric matrices with eigenvalues
$\lambda$ is as smooth as $f$ and symmetric of homogeneous degree
$m$, if $f$ is so. We denote by $(\dot{F}^{ij})$ the matrix of the
first partial derivatives of $F$ with respect to the components of
its arguments:
\begin{eqnarray*}\frac {\partial }{\partial s}F(A+sB)\big|_{s=0}
=\dot{F}^{ij}(A)B_{ij},\end{eqnarray*} 
where $A$ and $B$ are any symmetric matrices.
Similarly for the second partial derivatives of $F$, we write
\begin{eqnarray*}\frac {\partial ^2}{\partial s^2}F(A+sB)\big|_{s=0}
=\ddot{F}^{ij, kl}(A)B_{ij}B_{kl}.\end{eqnarray*} 
We also use the notation
\begin{eqnarray*}\dot f_i(\lambda )=\frac{\partial f}{\partial \lambda
_i}(\lambda), \mbox{ and } \ddot f_{ij}(\lambda )=\frac{\partial^2
f}{\partial \lambda _i\partial \lambda _j}(\lambda). \end{eqnarray*}

Recall that the homogeneity of $F=F(h_{ij})$ implies the following
\begin{equation}
\dot{F}^{ij}h_{ij}=mF \quad \mbox{and}\quad
\ddot{F}^{ij,rs}h_{ij}h_{rs}=(m-1)\dot{F}^{rs}h_{rs}.
\end{equation}

The following proposition is well known (see e.g. \cite{a4,g})\\[3mm]
{\bf Proposition 2.1}~ {\it Let $f$ and $F$ be as above. If $f$ is
$C^2$ and symmetric, then at any diagonal matrix $A$ with distinct
eigenvalues, the second-order derivative of $F$ in direction $B$ is
given by
\begin{eqnarray*}{\ddot F}(B,B)=\sum _{k,l}\ddot
f_{kl}B_{kk}B_{ll}+2\sum_{k<l}\frac{\dot f_k-\dot f_l}{\lambda
_k-\lambda _l}B_{kl}^2.
\end{eqnarray*}}\\[3mm]
\indent The following corollary follows immediately,\\[3mm]
{\bf Corollary 2.2}~ {\it If $f$ is convex (concave) at $\lambda
(A)$, then $F$ is convex (concave) at $A$. Moreover  $f$ is convex
(concave) if and only if
\begin{eqnarray*}\frac{\dot f_i-\dot f_j}{\lambda_{i}-\lambda_{j}} 
\geq (\leq)~ 0 ~~for ~all~ i \neq j.\end{eqnarray*}}
\\[3mm]
\indent Let $\mathcal{F}(A)=\mathfrak{f}(\lambda (A))$ be another
homogeneous function defined in $\R^n$. The first part of the next
lemma is in fact in \cite{lyw}, where only the symmetric and
homogeneous degree one function is considered.\\[3mm]
{\bf Lemma 2.3}~ {\it Assume $F$ and $\mathcal{F}$ are elliptic and
of homogeneous degree $m$, and the eigenvalues $\lambda (A)$ of $A$
are non-negative. If $f$ is convex (concave), and $\mathfrak{f}$
concave (convex), then
$$m(\mathcal{F}\dot{F}^{ij}h_{ik}h_{j}^{k}
-F\dot{\mathcal{F}}^{ij}h_{ik}h_{j}^{k})\geq (\le )0,$$ and
$$m\sum _j(F\dot{\mathfrak{f}}_j-\mathcal{F}\dot{f}_j)\geq (\le
)0.$$}
\begin{proof} For the first inequality, using the homogeneity of $F$ and 
$\mathcal{F}$, we compute as in \cite{lyw}
\begin{eqnarray*}&&m(\mathcal{F}\dot{F}^{ij}h_{ik}h_{j}^{k}
-F\dot{\mathcal{F}}^{ij}h_{ik}h_{j}^{k})=\sum _{i,j}(
\dot{\mathfrak{f}}_j\dot f_i\lambda _j\lambda _i^2-
\dot{\mathfrak{f}}_i\dot f_j\lambda _j\lambda _i^2)\\
&=&\frac 12\sum _{i\neq j}\lambda _i\lambda _j(\lambda _i-\lambda
_j)^2\left[\dot
{\mathfrak{f}}_j\left(\frac{\dot{f}_i-\dot{f}_j}{\lambda_i-\lambda
_j}\right)-\dot{f}_j\left(\frac{\dot {\mathfrak{f}}_i-\dot
{\mathfrak{f}}_j}{\lambda _i-\lambda _j}\right)\right].
\end{eqnarray*}
The lemma now follows by using Corollary 2.2. For the second
inequality, we similarly have
$$m\sum _j(F\dot{\mathfrak{f}}_j-\mathcal{F}\dot{f}_j)=
\frac 12\sum_{i,j}\left(\dot{\mathfrak{f}}_j(\lambda_i-\lambda _j)
(\dot{f}_i-\dot{f}_j)-\dot{f}_i(\lambda_i-\lambda_j)
(\dot{\mathfrak{f}}_i-\dot{\mathfrak{f}}_j)\right),$$
the required inequality follows. 
\end{proof}

Since $N$ is of constant curvature $c$, we have the Codazzi equation
$$\nabla _kh_{ij}=\nabla_jh_{ik}.$$

The Codazzi's equation implies the Ricci identity
\begin{eqnarray}
\nabla _i\nabla _jh_{kl}&=&\nabla _k\nabla
_lh_{ij}+h_{il}A_{kj}^2-h_{kl}A_{ij}^2+h_{ij}A_{kl}^2-h_{kj}A_{il}^2\nonumber\\
&&-c(g_{il}h_{kj}-g_{kl}h_{ij}+g_{ij}h_{kl}-g_{kj}h_{il}),
\end{eqnarray}
where $A_{ij}^s=h_i^lh_l^k\cdots h_{kj}$ ($s$ factors).

\section{Computations on Curvature Functions}
\renewcommand{\thesection}{\arabic{section}}
\renewcommand{\theequation}{\thesection.\arabic{equation}}
\setcounter{equation}{0}
\setcounter{theorem}{0}

Let $\{\frac{\partial}{\partial x^i}\}$ be the natural frame field
on $M^n$. Denote by $\bar{\nabla}$ the covariant derivative of
$N^{n+1}(c)$. The following is well-known
\begin{eqnarray}
<\bar{\nabla}_X\partial_{\rho}, Y>=\bar{\nabla}^2\rho (X,
Y)=\left\{\begin{array}{ll}0& \mbox{ if } X=\partial _{\rho}\\
\frac{ch_c(\rho)}{sh_c(\rho)}<X, Y>& \mbox{ if } <X, \partial
_{\rho}>=0\end{array}\right..
\end{eqnarray}

By (3.1) and Codazzi equation, we have the following lemma
\cite{cm}\\[3mm]
{\bf Lemma 3.1}~ {\it  The second order derivative of $\mathcal{Z}$
is given by
\begin{eqnarray*}
\nabla_i\nabla_j \mathcal{Z}=-ch_c(\rho)h_{ij}-sh_c(\rho)<\partial
_{\rho}^{\top},\nabla h_{ij}>-\mathcal{Z}A_{ij}^2,
\end{eqnarray*}
where $\partial _{\rho}^{\top}$ is the component of
$\partial_{\rho}$ tangent to $M$.}\\[3mm]
\indent
 We differentiate the equation (1.6) to get
$$\nabla_jF= \dot{F}^{kl} \nabla_jh_{kl}= -\tau \nabla_j\mathcal{Z}.$$

Taking derivative of the above equation again in a normal coordinate
system with the help of Lemma 3.1, we have
\begin{eqnarray*}
&&\ddot{F}^{kl,rs}\nabla _ih_{rs}\nabla_jh_{kl}
+\dot{F}^{kl}\nabla_i\nabla_jh_{kl} = \nabla_i\nabla_jF
=-\tau \nabla_i\nabla_j\mathcal{Z}\\
  &= & -\tau \left(-ch_c(\rho)h_{ij}-sh_c(\rho)<\partial_{\rho}^{\top},
\nabla h_{ij}>-\mathcal{Z}A_{ij}^2\right)
\\
& = & \tau ch_c(\rho) h_{ij}-FA_{ij}^2+\tau
sh_c(\rho)<\partial_{\rho}, \frac{\partial}{\partial
x^l}>\nabla^lh_{ij},
\end{eqnarray*}
which implies
\begin{equation}
\dot{F}^{kl}\nabla_i\nabla_jh_{kl}=\tau ch_c(\rho)
h_{ij}-FA^2_{ij}+\tau sh_c(\rho)<\partial_{\rho},
\frac{\partial}{\partial x^l}>\nabla^lh_{ij}-\ddot{F}^{kl,rs}\nabla
_ih_{rs}\nabla _jh_{kl}.
\end{equation}

 Using the Euler relation (2.1), we have by (3.2)
\begin{eqnarray}
\dot{F}^{kl}\nabla _k\nabla _lF
&=&\dot{F}^{kl}\nabla_k(\dot{F}^{ij}\nabla _lh_{ij})\nonumber\\
& = & \dot{F}^{kl}\ddot{F}^{ij,rs}\nabla _lh_{ij}\nabla_kh_{rs}
+\dot{F}^{ij}\dot{F}^{kl}\nabla _k\nabla_lh_{ij}\nonumber \\
& = & \dot{F}^{kl}\ddot{F}^{ij,rs}\nabla _lh_{ij}\nabla
_kh_{rs}+\dot{F}^{kl}\left(
 \tau ch_c(\rho)h_{kl}-FA^2_{kl} \right.\nonumber\\
&&\left.+ \tau sh_c(\rho)<\partial_{\rho},\frac{\partial}{\partial
x^i}>\nabla^ih_{kl}
-\ddot{F}^{ij,rs}\nabla _kh_{rs}\nabla_lh_{ij}\right)\nonumber\\
& = & m\tau ch_c(\rho)F+\tau sh_c(\rho)<\partial_{\rho},
\frac{\partial}{\partial x^l}>\nabla^lF -F\dot{F}^{kl}A_{kl}^2.
\end{eqnarray}
Here $m$ is the degree of $F$.

 For any other curvature function
$\mathcal{F}$ of homogeneous degree $m$, we compute similarly
\begin{eqnarray*}
\nabla_k\nabla_l\mathcal{F}=\ddot{\mathcal{F}}^{ij,rs}\nabla
_kh_{rs}\nabla
_lh_{ij}+\dot{\mathcal{F}}^{ij}\nabla_k\nabla_lh_{ij}.
\end{eqnarray*}

Using the Ricci identity (2.2) and inserting (3.2) into the above
equation, we obtain
\begin{eqnarray}
\dot{F}^{kl}\nabla_k\nabla_l\mathcal{F}&=&\dot{F}^{kl}
\ddot{\mathcal{F}}^{ij,rs}\nabla_kh_{rs}\nabla
_lh_{ij}+\dot{\mathcal{F}}^{ij}\dot{F}^{kl}\nabla_k\nabla_lh_{ij}\nonumber\\
&=&\dot{F}^{kl}\ddot{\mathcal{F}}^{ij,rs}\nabla _kh_{rs}\nabla
_lh_{ij}\nonumber\\
&&+\dot{\mathcal{F}}^{ij}\left[\dot{F}^{kl}\nabla_i\nabla_jh_{kl}-\dot{F}^{kl}
(h_{il}A_{kj}^2-h_{kl}A_{ij}^2+h_{ij}A_{kl}^2-h_{kj}A_{il}^2)\right]\nonumber\\
&&\;\;\;\;\;\;\;\;\;+c\dot{\mathcal{F}}^{ij}\dot{F}^{kl}(g_{il}h_{kj}
-g_{kl}h_{ij}+g_{ij}h_{kl}-g_{kj}h_{il})\nonumber\\[2mm]
&=&\dot{F}^{kl}\ddot{\mathcal{F}}^{ij,rs}\nabla _kh_{rs}\nabla
_lh_{ij}+\dot{\mathcal{F}}^{ij}(mFA_{ij}^2-h_{ij}\dot{F}^{kl}A_{kl}^2)
-cm(\mathcal{F}\dot F^{kk}-F\dot {\mathcal{F}}^{kk})\nonumber\\
&&+\dot{\mathcal{F}}^{ij}\left[\tau ch_c(\rho)h_{ij}-FA_{ij}^2+\tau
sh_c(\rho)<\partial_{\rho}, \frac{\partial}{\partial
x^l}>\nabla^lh_{ij}-\ddot{F}^{kl,rs}\nabla
_ih_{rs}\nabla _jh_{kl}\right]\nonumber\\
&=&m\tau ch_c(\rho)\mathcal{F}+\tau sh_c(\rho)<\partial_{\rho},
\frac{\partial}{\partial x^l}>\nabla^l
\mathcal{F}+(m-1)F\dot{\mathcal{F}}^{ij}A_{ij}^2-m\mathcal{F}
\dot{F}^{ij}A_{ij}^2\nonumber\\
&&+(\dot{F}^{ij}\ddot{\mathcal{F}}^{kl,rs}-
\dot{\mathcal{F}}^{ij}\ddot{F}^{kl,rs})\nabla
_ih_{rs}\nabla _jh_{kl}-cm(\mathcal{F}\dot F^{kk}-F\dot
{\mathcal{F}}^{kk}).
\end{eqnarray}

Direct computation gives
\begin{eqnarray*}
 \dot{F}^{kl}\nabla _k\nabla _l \mbox{\Large{$($}}
\frac {F}{\mathcal{F}}\mbox{\Large{$)$}} &=& \frac
{1}{\mathcal{F}}\dot{F}^{kl}\nabla _k\nabla _lF- \frac
{F}{\mathcal{F}^{2}}\dot{F}^{kl}\nabla _k\nabla _l\mathcal{F}\\
&&-\frac  2{\mathcal{F}^{2}}\dot{F}^{kl}\nabla _kF\nabla
_l\mathcal{F} +\frac {2F}{\mathcal{F}^{3}}\dot{F}^{kl}\nabla
_k\mathcal{F}\nabla_l\mathcal{F},
\end{eqnarray*}
which implies by (3.3) and (3.4)
\begin{eqnarray}
\dot{F}^{kl}\nabla _k\nabla _l\left(\frac {F}{\mathcal{F}}\right)&
=& \frac {1}{\mathcal{F}}\left(m\tau ch_c(\rho)F+\tau
sh_c(\rho)<\partial_{\rho}, \frac{\partial}{\partial
x^l}>\nabla^lF -F\dot{F}^{kl}A_{kl}^2\right)\nonumber\\
&&-\frac {F}{\mathcal{F}^{2}}\left[m\tau ch_c(\rho)\mathcal{F}+\tau
sh_c(\rho) <\partial_{\rho}, \frac{\partial}{\partial x^l}>\nabla^l
\mathcal{F}\right.\nonumber\\
&&~~~~~~~~~~~+(m-1)F\dot{\mathcal{F}}^{ij}A_{ij}^2-m\mathcal{F}
\dot{F}^{ij}A_{ij}^2\nonumber\\
&&+\left.(\dot{F}^{ij}\ddot{\mathcal{F}}^{kl,rs}-
\dot{\mathcal{F}}^{ij}\ddot{F}^{kl,rs})\nabla
_ih_{rs}\nabla _jh_{kl}\right.\left.-cm(\mathcal{F}\dot F^{kk}-F\dot
{\mathcal{F}}^{kk})\right]\nonumber\\
&&-\frac 2{\mathcal{F}^{2}}\dot{F}^{kl}\nabla _kF\nabla
_l\mathcal{F} +\frac {2F}{\mathcal{F}^{3}}\dot{F}^{kl}\nabla
_k\mathcal{F}\nabla_l\mathcal{F}.
\end{eqnarray}

Note that the last two terms in  (3.5) are equal to $-\frac
    2{\mathcal{F}}\dot{F}^{kl}\nabla _k\mathcal{F}
\nabla _l\mbox{\Large{$($}}\frac
    {F}{\mathcal{F}}\mbox{\Large{$)$}}$,
and so we at last arrive at the lemma\\[3mm]
{\bf Lemma 3.2}~ {\it Let $F$ and $\mathcal{F}$ be two nonzero
curvature functions on $M$, which are homogeneous of degree $m$. If
$F$ satisfies (1.6), then the following holds
\begin{eqnarray}
 \dot{F}^{kl}\nabla _k\nabla _l\left(\frac {F}{\mathcal{F}}\right)
 & = &   \tau sh_c(\rho)<\partial_{\rho}, \frac{\partial}{\partial x^l}
>\nabla^l\left(\frac {F}{\mathcal{F}}\right)-\frac
    2{\mathcal{F}}\dot{F}^{kl}\nabla _k\mathcal{F}\nabla _l\left(\frac
    {F}{\mathcal{F}}\right)\nonumber\\
   &&+\frac{F}{\mathcal{F}^2}(\dot{\mathcal{F}}^{ij}\ddot{F}^{kl,rs}
-\dot{F}^{ij}\ddot{\mathcal{F}}^{kl,rs})\nabla_ih_{rs}\nabla_jh_{kl}\\
&&+(m-1)\frac{F}{\mathcal{F}^2}(\mathcal{F}\dot{F}^{ij}A_{ij}^2
-F\dot{\mathcal{F}}^{ij}A_{ij}^2)\\
&&-\frac {cmF}{\mathcal{F}^2}(F\dot
{\mathcal{F}}^{kk}-\mathcal{F}\dot F^{kk}).
\end{eqnarray}}

\section{Proof of the Main Theorem}
\renewcommand{\thesection}{\arabic{section}}
\renewcommand{\theequation}{\thesection.\arabic{equation}}
\setcounter{equation}{0}
\setcounter{theorem}{0}

Firstly we consider the case $m\ge 1$. It is clear, in this case,
$F$ is positive restricting to $M$. We only prove the theorem for
$f$ concave. It is similar for $f$ convex. Taking an elliptic convex
curvature function $\mathcal{F}(A)=\mathfrak{f}(\lambda (A))$ of
homogeneous degree $m$, the homogeneity implies that $\mathcal{F}$
is also positive as $F$. By Corollary 2.2, Lemma 2.3 and the
homogeneity of $F$ and $\mathcal{F}$, we see that (3.6), (3.7) and
(3.8) are non-positive since $c\le 0$. Then applying the strong
maximum principle to $\frac {F}{\mathcal{F}}$ in Lemma 3.2 yields
$\frac{F}{\mathcal{F}}=c_1$, a positive constant, on $M$. Therefore
by assumption, either $F$ and $\mathcal{F}$ are constant restricting
to $M$ or
$$0\ge \ddot{F}^{ij,kl}\eta_{ij}\eta_{kl}=c_1\ddot{\mathcal{F}}^{ij,kl}
\eta_{ij}\eta_{kl}\ge 0,$$
for any real symmetric matrix $\eta$. Especially,
$$\dot{F}^{ij}\ddot{F}^{kl,rs}\nabla
_ih_{rs}\nabla _jh_{kl}=0,$$ and
\begin{eqnarray}
\ddot{F}^{kl,rs}h_{rs}h_{kl}=0.
\end{eqnarray}

If $\nabla _ih_{kl}=0$ for any $i,j,k=1, 2, \cdots, n$,  then the
mean curvature $H$ is constant and we have done. Otherwise by (2.1),
(4.1) implies $m(m-1)F=0$, and so $m=1$. Since $H$ is concave as
well as convex, taking $\mathcal{F}=H$, we have $F=c_1H$. The
theorem  now follows from  Proposition 4.1 below.

Secondly, we consider the case $m<0$. Again, the homogeneity and
ellipticity of $F$ imply that $F<0$ since $m<0$. As in the first
case, we only consider the case for $f$ concave, and it is similar
for $f$ convex. As before we take an elliptic and convex curvature
function $\mathcal{F}$ which is homogeneous of degree $m$. Then
terms in (3.6)-(3.8) are non-negative. Applying again the strong
maximum principle to $\frac{F}{\mathcal{F}}$ in Lemma 3.2 yields
$\frac{F}{\mathcal{F}}=c_2$, a positive constant. Similar discussion
also as in the first case, we see the only possible case is $\nabla
_ih_{kl}=0$ for any $i,j,k=1, 2, \cdots, n$, and therefore the mean
curvature $H$ is a constant on $M$,  which implies $X(M)$ is again a
sphere.

Lastly, we consider the case $n=2$. For any symmetric function
$\mathcal{F}$ of the principal curvatures, which is homogeneous of
degree zero, we compute as in (3.4) to obtain
\begin{eqnarray}
\dot{F}^{kl}\nabla_k\nabla_l\mathcal{F}&=&\dot{F}^{kl}
\ddot{\mathcal{F}}^{ij,rs}\nabla_kh_{rs}\nabla
_lh_{ij}+\dot{\mathcal{F}}^{ij}\dot{F}^{kl}\nabla_k\nabla_lh_{ij}\nonumber\\
&=&\tau sh_c(\rho)<\partial_{\rho}, \frac{\partial}{\partial
x^l}>\nabla^l
\mathcal{F}+(m-1)F\dot{\mathcal{F}}^{ij}A_{ij}^2\nonumber\\
&&+(\dot{F}^{ij}\ddot{\mathcal{F}}^{kl,rs}-\dot{\mathcal{F}}^{ij}
\ddot{F}^{kl,rs})\nabla_ih_{rs}\nabla _jh_{kl}+cmF\dot {\mathcal{F}}^{kk}.
\end{eqnarray}
We now compute the second-order derivatives in terms of Proposition
2.1. Since $n=2$, it's not difficult to check as in \cite{a5} that,
the terms in (4.2) containing the second-order derivatives of $F$
and $\mathcal{F}$ are given by in a frame diagonalizing the second
fundamental form
\begin{eqnarray*}
Q&=&(\dot{F}^{ij}\ddot{\mathcal{F}}^{kl,rs}-\dot{\mathcal{F}}^{ij}
\ddot{F}^{kl,rs})\nabla_ih_{rs}\nabla _jh_{kl}\nonumber\\
&=&\;\;\;\;(\dot{f}_1\ddot{\mathfrak{f}}_{11}-\dot{\mathfrak{f}}_1
\ddot{f}_{11})(\nabla_1h_{11})^2
+(\dot{f}_1\ddot{\mathfrak{f}}_{22}-\dot{\mathfrak{f}}_1\ddot{f}_{22})
(\nabla_1h_{22})^2\nonumber\\[2mm]
&&+(\dot{f}_2\ddot{\mathfrak{f}}_{11}-\dot{\mathfrak{f}}_2\ddot{f}_{11})
(\nabla_2h_{11})^2
+(\dot{f}_2\ddot{\mathfrak{f}}_{22}-\dot{\mathfrak{f}}_2\ddot{f}_{22})
(\nabla_2h_{22})^2\nonumber\\[2mm]
&&+2(\dot{f}_1\ddot{\mathfrak{f}}_{12}-\dot{\mathfrak{f}}_1\ddot{f}_{12})
\nabla_1h_{11}\nabla_1h_{22}
+2(\dot{f}_2\ddot{\mathfrak{f}}_{12}-\dot{\mathfrak{f}}_2\ddot{f}_{12})
\nabla_2h_{11}\nabla_2h_{22}\nonumber\\[2mm]
&&+2\frac{\dot{f}_1\dot{\mathfrak{f}}_2-\dot{f}_2\dot{\mathfrak{f}}_1}
{\lambda_2-\lambda_1}(\nabla_1h_{12})^2
+2\frac{\dot{f}_1\dot{\mathfrak{f}}_2-\dot{f}_2\dot{\mathfrak{f}}_1}
{\lambda_2-\lambda_1}(\nabla_2h_{12})^2.
\end{eqnarray*}
As in \cite{a5} again, we can work at a maximum point of
$\mathcal{F}$. Then using the gradient conditions and the
homogeneity of $F$, we have
\begin{eqnarray}
Q&=&-mF\dot{\mathfrak{f}}_1\left(\frac{m-1}{\lambda_2^2}+\frac{2}
{\lambda_1(\lambda_2-\lambda_1)}\right)
(\nabla_1h_{22})^2\nonumber\\
&&-mF\dot{\mathfrak{f}}_2\left(\frac{m-1}{\lambda_1^2}-\frac{2}
{\lambda_2(\lambda_2-\lambda_1)}\right)
(\nabla_2h_{11})^2.
\end{eqnarray}
By Euler identity we also compute
$$F\dot{\mathcal{F}}^{ij}A_{ij}^2=F\dot{\mathfrak{f}}_2\lambda_2
(\lambda_2-\lambda_1).$$
Similarly for the last term in (4.2)
$$cmF\dot {\mathcal{F}}^{kk}=-cmF\dot{\mathfrak{f}}_2
\frac{\lambda_2-\lambda_1}{\lambda_1}.$$
Putting these formulae into (4.2), we have at a maximum point of
$\mathcal{F}$
\begin{eqnarray}
\dot{F}^{kl}\nabla_k\nabla_l\mathcal{F}&=&\tau
sh_c(\rho)<\partial_{\rho}, \frac{\partial}{\partial x^l}>\nabla^l
\mathcal{F}\nonumber\\
&&+(m-1)F\dot{\mathfrak{f}}_2\lambda_2(\lambda_2-\lambda_1)
-cmF\dot{\mathfrak{f}}_2\frac{\lambda_2-\lambda_1}{\lambda_1}\\
&&-mF\dot{\mathfrak{f}}_1\left(\frac{m-1}{\lambda_2^2}
+\frac{2}{\lambda_1(\lambda_2-\lambda_1)}\right)
(\nabla_1h_{22})^2\\
&&-mF\dot{\mathfrak{f}}_2\left(\frac{m-1}{\lambda_1^2}
-\frac{2}{\lambda_2(\lambda_2-\lambda_1)}\right)(\nabla_2h_{11})^2.
\end{eqnarray}

Now take
$\mathcal{F}=\frac{2|A|^2-H^2}{H^2}=\frac{(\lambda_1-\lambda_2)^2}
{(\lambda_1+\lambda_2)^2}$,
we have
$\dot{\mathfrak{f}}_1=\frac{4\lambda_2(\lambda_1-\lambda_2)}
{(\lambda_1+\lambda_2)^3}$,
and
$\dot{\mathfrak{f}}_2=-\frac{4\lambda_1(\lambda_1-\lambda_2)}
{(\lambda_1+\lambda_2)^3}$.
We assume $\lambda _1\le \lambda_2$ at the maximum point of
$\mathcal{F}$. Then when $m\ge 1$ and $c\le 0$, (4.4) is
non-negative. If $m=1$, (4.5) and (4.6) are all non-negative, we
have immediately by strong maximum principle $\mathcal{F}$ is a
constant, and therefore $X(M)$ is an umbilical sphere.

If $m>1$, in order to apply the maximum principle, we require
$\frac{m-1}{\lambda_2^2}+\frac{2}{\lambda_1(\lambda_2-\lambda_1)}$
is non-negative, and
$\frac{m-1}{\lambda_1^2}-\frac{2}{\lambda_2(\lambda_2-\lambda_1)}$
is non-positive. Thus the pinching ratio
$r=\frac{\lambda_2}{\lambda_1}$ must satisfy the conditions
\begin{eqnarray}2r^2+(m-1)r-(m-1)\ge 0,
\end{eqnarray} and
\begin{eqnarray}(m-1)r^2-(m-1)r-2\le 0.
\end{eqnarray}
 The first is always true since $r\ge
1$, and the second is true if and only if
$$r\le\frac12 \left(1+\sqrt{1+\frac {8}{m-1}}\right).$$
Then when $r$ satisfies the above inequality, by maximum principle,
$\mathcal{F}$ is a constant and therefore $X(M)$ is an umbilical
sphere.

If $m<0$, since (4.5) and (4.6) are non-negative,  we also require
(4.7) and (4.8) hold since $F<0$. It is easy to check that when
$-7\le m<0$, (4.7) and (4.8) are always satisfied. For $m<-7$, (4.8)
is always true, and (4.7) is true if and only if
$$r\le\frac14 \left((1-m)-\sqrt{(1-m)^2-8(1-m)}\right)=
\frac{2}{1+\sqrt{1-\frac{8}{1-m}}}.$$
Therefore when $-7\le m<0$ or when $m<-7$ and $r\le
2/\left(1+\sqrt{1-\frac{8}{1-m}}\right)$, the maximum principle
implies that $X(M)$ is an umbilical sphere.
 $\hfill \Box$\\

When $c=0$, the following Proposition 4.1 is essentially a result of
Huisken \cite{h3}, for it differs from his result only by a constant
$\tau$. For completeness we give the proof.\\[3mm]
{\bf Proposition 4.1}~ {\it If $X: M^n\rightarrow N^{n+1}(c\le 0)$
is compact, connected, with non negative mean curvature and
satisfies $H+\tau \mathcal{Z}=0$ for some positive constant $\tau $
depending only on $n$, then $X(M)$ is an umbilical sphere.}
\begin{proof}
By (3.2) with $F=H$, we have
\begin{equation}\triangle H=\tau ch_c(\rho)
H -H|A|^2+\tau sh_c(\rho)<\partial_{\rho}, \frac{\partial}{\partial
x^l}>\nabla ^lH,
\end{equation} which implies that
$H>0$ by strong maximum principle, and from Ricci identity (2.2)
\begin{equation}
\triangle |A|^2=2|\nabla A|^2+2\tau ch_c(\rho)|A|^2 -2|A|^4+\tau
sh_c(\rho)<\partial_{\rho}, \frac{\partial}{\partial x^l}>\nabla
^l|A|^2-2c(n|A|^2-H^2).
\end{equation}
Using (4.9) and (4.10), by similar calculation as in section 3, we
have
\begin{eqnarray}
\triangle \mbox{\Large{$($}} \frac {|A|^2}{|H|^2}\mbox{\Large{$)$}}
&=& \tau sh_c(\rho)<\partial_{\rho}, \frac{\partial}{\partial
x^l}>\nabla ^l\mbox{\Large{$($}}\frac {|A|^2}{H^2}\mbox{\Large{$)$}}
+\frac 2{H^4}\,|h_{ij}\nabla _lH-\nabla
_lh_{ij}H|^2\nonumber\\
&&-\frac 1{H^2}<\nabla H^2, \nabla \mbox{\Large{$($}}\frac
{|A|^2}{H^2}\mbox{\Large{$)$}}>-\frac {2c}{H^2}(n|A|^2-H^2).
\end{eqnarray}
Since $M$ is compact, the strong maximum principle implies that
\begin{equation}
\frac {|A|^2}{H^2}=\mbox {constant} \quad \mbox {and} \quad
|H\nabla_ih_{kl}-\nabla_iHh_{kl}|^2\equiv 0.
\end{equation}
Then if $c=0$, Huisken's theorem implies that $X(M)$ is a sphere. If
$c<0$, we have by (4.11) and (4.12), $n|A|^2-H^2=0$. It follows
immediately that $X(M)$ is an umbilical sphere. 
\end{proof}


\begin{thebibliography}{99}
{
\bibitem{a1} Andrews B., \em  Contraction of convex hypersurfaces in
Euclidean space, \em  Calc.\ Var.\ Partial Differential Equations.
{\bf 2} (1994), no.\ 2, 151--171.
\bibitem{a2} Andrews B., {\em Evolving convex curves}.
Calc. Var. Partial Differential Equations. {\bf 7} (1998), no.\ 3, 315--371.
\bibitem{a3} Andrews B., {\em Classification of limiting shapes for
isotropic curve flows}.  J.\ Amer.\ Math.\ Soc. {\bf 16} (2003), no.\ 2, 
 443--459.
\bibitem{a4} Andrews B., {\em Pinching estimates and motion of
hypersurfaces by curvature functions}. J.\ reine angew.\ Math.\
{\bf 608} (2007), 17--33.
\bibitem{a5} Andrews B., {\em Moving surfaces by non-concave curvature
functions}. arXiv: math.DG/0402273, 2004.
\bibitem{cgm} Carreras F., Gimenez F. and Miquel V.,
{\em Immersions of compact Riemannian manifolds into a ball of a complex
space form}. Math.\ Z.\ {\bf 225} (1997), 103--113.
\bibitem{cm} Cabezas-Rivas E. and Miquel V., {\em Volume preserving mean
curvature flow in the hyperbolic space}.
Indiana Univ.\ Math.\ J.\ {\bf 56} (2007), no.\ 5, 2061--2086.
\bibitem{ct} Chow  B. and Tsai  H., {\em Geometric expansion of convex
plane curves}. J.\ Differential Geom. {\bf 44} (1996), 312--330.
\bibitem{eh} Ecker K. and Huisken G., {\em Immersed hypersurfaces
with constant Weingarten curvature}. Math.\ Ann.\ {\bf 283} (1989), 329--332.
\bibitem{h1} Huisken G., {\em Flow by mean curvature of convex surfaces
into spheres}. J.\ Differential Geom. {\bf 20} (1984), no.\ 1, 237--266.
\bibitem{h3} Huisken G., {\em Asymptotic behavior for singularities
of the mean curvature flow}. J.\ Differential Geom.\ {\bf 31} (1990),
no.\ 1, 285--299.
\bibitem{g} Gerhardt C., {\em Closed Weingarten hypersurfaces in
Riemannian manifolds}. J.\ Differential Geom.{\bf 43} (1996), no.\ 3, 
612--641.
\bibitem{lyw} Li G., Yu L.\ and Wu C., {\em Curvature flow with a general
forcing term in Euclidean spaces}. J.\ Math.\ Anal.\ Appl.\
 {\bf 353} (2009), 508--520.
\bibitem{m3} McCoy J., {\em Mixed volume preserving curvature flows}.
 Calc.\ Var.\ Partial Differential Equations {\bf 24} (2005), 131--154.
\bibitem{r} Ros A., {\em Compact hypersurfaces with constant scalar curvature
and a congruence theorem}, J.\ Differential Geom. {\bf 27} (1988), 215--220.
\bibitem{s} Smoczyk K., {\em Self-shrinkers of the mean curvature flow in
arbitrary codimension}. Int.\ Math.\ Res.\ Not.\ {\bf 48} (2005), 2983--3004.
\bibitem{u} Urbas J., {\em Convex curves moving homothetically by negative 
powers of their curvature}. Asian J.\  Math.\ {\bf 3} (1999), 635--656.
\bibitem{z} Zhu X., {\em Lectures on Mean Curvature Flow}. 
Studies on Advanced Mathematics, American Mathematical Society, 
International Press, {\bf 32} (2002).}
\end{thebibliography}
 \end{document}